\documentclass{amsart}
\usepackage{amssymb}
\usepackage{graphicx} 
\usepackage{eepic}

\newtheorem{theorem}{Theorem}[section]

\theoremstyle{definition}

\newtheorem{corollary}[theorem]{Corollary}
\newtheorem{proposition}[theorem]{Proposition}
\newtheorem{example}[theorem]{Example}

\theoremstyle{remark}

\numberwithin{equation}{section}

\newcommand{\Z}{ {\mathbb Z} }

\newcommand{\gen}[1]{\left\langle #1 \right\rangle}

\begin{document}
\bibliographystyle{plain}
\title[]{Shellable Complexes from Multicomplexes}
\author{Jonathan Browder}
\date{24 December 2008}

\begin{abstract}
 Suppose a group $G$ acts properly on a simplicial complex $\Gamma$. Let $l$ be the number of $G$-invariant vertices and $p_1, p_2, \ldots p_m$ be the sizes of the $G$-orbits having size greater than 1. Then $\Gamma$ must be a subcomplex of $\Lambda = \Delta^{l-1}* \partial \Delta^{p_1-1}*\ldots * \partial \Delta^{p_m-1}$. A result of Novik gives necessary conditions on the face numbers of Cohen-Macaulay subcomplexes of $\Lambda$. We show that these conditions are also sufficient, and thus provide a complete characterization of the face numbers of these complexes.

\end{abstract}

\maketitle

\section{Introduction}
\label{intro}

One of the central problems in geometric combinatorics is that of characterizing the face numbers of various classes of simplicial complexes. The Kruskal-Katona theorem \cite{MR0154827, MR0290982} characterized the $f$-vectors of all simplicial complexes, while a result of Stanley characterized the face numbers of all Cohen-Macaulay complexes \cite{MR0572989}. One fruitful line of inquiry since then has been in determining additional conditions on the face numbers of complexes with certain types of symmetry.

 In particular, let $\Gamma$ be a simplicial complex on $n$ vertices, and suppose $G$ is a group which acts on $\Gamma$. We say the action of $G$ is \emph{proper} if whenever $F$ is a face of $\Gamma$ and $gF = F$ for some $g \in G$, then $gv = v$ for each vertex $v \in F$, i.e., whenever an element of $G$ fixes a face of $\Gamma$ it fixes that face pointwise. Let $V'$ be the set of $G$-invariant vertices of $\Gamma$ and let $V_1, V_2, \ldots , V_m$ be the $G$-orbits on the vertex set of $\Gamma$ with size greater than $1$. If the action of $G$ is proper, no face of $\Gamma$ can contain any $V_i$, so $\Gamma$ must be a subcomplex of $\Lambda(l;p_1,p_2,\ldots ,p_m) = \Delta^{l-1}* \partial \Delta^{p_1-1}*\ldots * \partial \Delta^{p_m-1}$, where $l = |V'|$, $p_i = |V_i|$, $\Delta^{k}$ is the $k$-simplex and $\partial \Delta^k$ is the boundary complex of $\Delta^k$. (Note also that as each face of $\Gamma$ must miss at least one element of each $V_i$, the dimension of $\Gamma$ is  at most $n-m-1$.)
 
 Let $S(a_1, a_2, \ldots , a_k)$ (for $0\leq a_i \leq \infty$) denote the set of all monomials $x_1^{c_1}x_2^{c_2}\cdots x_k^{c_k}$ with $c_i \leq a_i$. For short, we will write $S(\infty^r, a_{r+1}, \ldots , a_k)$ for $S(a_1, a_2, \ldots , a_k)$ when $a_i = \infty$ for $1 \leq i \leq r$.  A non-empty subset $M$ of $S(a_1, a_2, \ldots , a_k)$ is called a \emph{multicomplex} if it is closed under divisibility; that is, if whenever $\mu|\mu'$ and $\mu' \in M$, then $\mu \in M$. For $M$ finite, let $\text{deg}(M) = \text{max}\{ \text{deg}(\mu) : \mu \in M \}$. The $F$-vector of a multicomplex $M$ is $F(M) = (F_0, F_1, F_2, \ldots)$ where $F_i$ is the number of elements in $M$ of total degree $i$.
 
Recall that the $h$-vector of a (d-1)-dimensional simplicial complex $\Gamma$ is $h(\Gamma) = (h_0, h_1, \ldots, h_d)$ defined by $\sum_{i=0}^d h_i x^i = \sum_{i=0}^d f_{i-1}x^i(1-x)^{d-i}$ where $f_i$ is the number of $i$-dimensional faces of $\Gamma$. In particular, the $h$-vector of $\Gamma$ completely determines the face numbers of $\Gamma$. 

 The following result is essentially due to Novik \cite{MR2122283}. (In fact Novik considered the case $p_i = p_j$ for all $i, j$, but with slight modifications her proof gives the general case, as we will address in section \ref{generalize}).
 
 \begin{theorem} \label{novik}
 Let $\Gamma$ be a (d-1)-dimensional Cohen-Macaulay complex on $n=l+\sum_{i=1}^mp_i$ vertices, where $p_1, p_2, \ldots, p_m \geq 2$, $m,l \geq 0$ are arbitrary integers. If $\Gamma$ is a subcomplex of $\Lambda(l;p_1,p_2,\ldots ,p_m)$, then there is a multicomplex $M \subseteq S(\infty^{n-d-m},(p_1-1),(p_2-1), \ldots , (p_m-1))$ such that the $h$-vector of $\Gamma$ is equal to the $F$-vector of $M$.
 \end{theorem}

 The goal of this paper is to show the converse to this theorem. In fact, we establish a slightly stronger result.
 
 \begin{theorem} \label{main} Let $l \geq 1$, $p_1, p_2, \ldots, p_m \geq 2$ be arbitrary integers. Let $n= l + \sum_{i=1}^m p_i$ and suppose $d \leq n-m$. If $M \subseteq S(\infty^{n-d-m},p_1-1, p_2-1, \ldots , p_m-1)$ is a multicomplex such that $\text{deg}(M) \leq d$, then there is a (d-1)-dimensional shellable subcomplex $\Gamma$ of $\Lambda(l;p_1,p_2,\ldots ,p_m)$ such that $h(\Gamma) = F(M)$.
 \end{theorem}
 
 Combined with Theorem \ref{novik}, this gives a generalization of a theorem of Stanley \cite{MR0572989}, which asserts that $h = (h_0, h_1, \ldots h_d)$ is the $h$-vector of a Cohen-Macaulay complex of dimension $d-1$ if and only if $h$ is the $F$ vector of some multicomplex $M \subseteq S(\infty^{n-d})$.
 
 \begin{corollary} \label{gen} Let $p_1, p_2, \ldots , p_m \geq 2$, $m,l \geq 0$ be arbitrary integers, $n = l+\sum_{i=1}^mp_i$, and $d \leq n-m$. Suppose $F = (F_0, F_1, \ldots , F_d)$. Then the following are equivalent:
 \begin{enumerate} 
 \item $F$ is the $h$-vector of a shellable subcomplex of $\Lambda(l;p_1,p_2,\ldots ,p_m)$. 
 \item $F$ is the $h$-vector of a Cohen-Macaulay subcomplex of $\Lambda(l;p_1,p_2,\ldots ,p_m)$.
 \item $F$ is the $F$-vector of a multicomplex in $S(\infty^{n-d-m},p_1-1,p_2-1, \ldots , p_m-1)$.
 \end{enumerate}
 \end{corollary}

Before moving on, we note that a different generalization of Stanley's theorem was obtained by Bj\"{o}rner, Frankl, and Stanley for balanced Cohen-Macaulay complexes \cite{MR905148}, which we state here for comparison. 

Partition the vertex set of a simplicial complex $\Gamma$ into $m$ disjoint subsets, $V = V_1 \cup V_2 \cup \ldots \cup V_m$, and let $a = (a_1, a_2, \ldots , a_m)$ be a positive integer vector. We say $\Gamma$ is $a$\emph{-balanced} if for each facet $\tau$ of $\Gamma$ and $1 \leq i \leq m$, $|\tau \cap V_i| = a_i$. Similarly, a multicomplex $M$ is \emph{colored} of type $a$ if its set of indeterminates can be partitioned into sets $X_1, X_2, \ldots , X_m$ such that for any monomial $m = x_1^{b_1}x^{b_2}_2 \ldots x^{b_k}_k \in M$ and $1 \leq i \leq m$, $\sum_{x_j \in X_i}b_j \leq a_i$ (that is, the part of $m$ supported in the variables in $X_i$ has degree less than or equal to $a_i$). For $b = (b_1, b_2, \ldots , b_m)$ and $0 \leq b \leq a$, define $f_b$ to be the number of faces of $\Gamma$ that contain exactly $b_i$ elements of $V_i$ for each $i$. The array $(f_b)_{0 \leq b \leq a}$ is the refined $f$-vector of $\Gamma$, the refined $h$-vector $(h_b)_{0 \leq b \leq a}$ is given by 
\begin{equation*}
h_b = \sum_{c \leq b}f_c \prod_{i=1}^s (-1)^{b_i-c_i}{a_i - c_i \choose b_i-c_i}.
\end{equation*}

Similary, the refined $F$-vector of $M$ is $(F_b)_{0\leq b \leq a}$ where $F_b$ is the number of monomials $m \in M$ such that part of $m$ supported in $X_i$ has degree $b_i$. 

\begin{theorem} \label{bfsgen} \cite{MR905148} Let $a = (a_1, a_2, \ldots , a_m)$ be a positive integer vector and suppose $F = (F_b)_{0\leq b \leq a}$ is an array of integers. Then the following are equivalent:
 \begin{enumerate} 
 \item $F$ is the refined $h$-vector of an $a$-balanced shellable complex. 
 \item $F$ is the refined $h$-vector of an $a$-balanced Cohen-Macaulay complex.
 \item $F$ is the refined $F$-vector of a multicomplex which is balanced of type $a$.
 \end{enumerate}
\end{theorem}

In particular, the proof of $(3) \Rightarrow (1)$  has a very similar structure to our proof of Theorem \ref{main}, and there seems to be a close relationship between the two results.

 \section{Idea of the Proof}
 \label{idea}
 For $\tau$ a face of some simplicial complex, denote by $\overline{\tau}$ the set of all subsets of $\tau$. Recall that a ($d-1$)-dimensional simplicial complex $\Gamma$ is \emph{shellable} if it is pure (i.e., all of its facets have dimension $d-1$) and there is an ordering of its facets $(\tau_1, \tau_2, \ldots \tau_r)$ such that for $1 < i \leq r$, the complex $\overline{\tau}_i \cap \left( \cup_{j <i} \overline{\tau}_j \right)$ is pure of dimension $d-2$. Such an ordering is then called a \emph{shelling} of $\Gamma$. For $L = (\tau_1, \tau_2, \ldots , \tau_r)$ any ordering of the facets of $\Gamma$, let $T_L(\tau_i)$ denote the set of facets of $\overline{\tau}_i \cap \left( \cup_{j <i} \overline{\tau}_j \right)$ (which will be some set of subsets of $\tau_i$ of size $d-1$ if $L$ is a shelling) for $i>1$, and set $T_L(\tau_1) = \emptyset$. We then have the following nice characterization of the $h$-vector of $\Gamma$:
 
 \begin{proposition} \cite{MR0301635} Let $(h_0, h_1, \ldots, h_d)$ be the $h$-vector of $\Gamma$. Then if $L$ is a shelling of $\Gamma$, $h_i = | \{ \tau_j: |T_L(\tau_j)| = i \}|$.
 
 \end{proposition}
 
Now, suppose $\Gamma$ is a simplicial complex with shelling $L$, and suppose $K$ is a subset of the set of facets of $\Gamma$. Let $L' = (\tau'_1, \tau'_2, \ldots , \tau'_{r'})$ be the ordering of $K$ inherited from $L$. Suppose that $T_{L'}(\tau) = T_L(\tau)$ for each $\tau \in K$. Then it follows immediately that $\Gamma'= \cup_{i=1}^{r'} \overline{\tau'_i}$ is a shellable subcomplex of $\Gamma$ with $h$-vector $(h'_1, h'_2, \ldots ,h'_d)$, where $h_i' = |\{ \tau \in K : T_L(\tau) = i \}|$.
 
 To prove Theorem \ref{main} we will construct a shelling $L$ of the ($d-1$)-skeleton, $skel_d(\Lambda)$, of $\Lambda(l;p_1,p_2,\ldots ,p_m)$, and show that for $(F_0, F_1, \ldots, F_d)$ the $F$-vector of some multicomplex $M$ in $S(\infty^{n-d-m},p_1-1, p_2-1, \ldots , p_m-1)$, there is a subsequence $L' = (\tau'_1, \ldots , \tau'_r)$ of $L$ such that each $T_{L'}(\tau'_i) = T_{L}(\tau'_i)$, and the number of $\tau'_i$ with $|T_L(\tau_i')| = j$ is $F_j$. We then have a shellable subcomplex of $skel_d(\Lambda)$ with $h$-vector equal to the $F$-vector of $M$.

 To do this, we will establish a bijection $\sigma$ between the set of facets of $skel_d(\Lambda)$ and the set $S^d$ of elements of $S(\infty^{n-d-m},p_1-1, p_2-1, \ldots , p_m-1)$ with degree less than or equal to $d$, with the property that $|T_L(\tau)| = \text{deg}(\sigma(\tau))$. For $M \subseteq S^d$ a multicomplex, let $L^M$ be the restriction of $L$ to $\sigma^{-1}(M)$. Then if $T_{L^M}(\tau) = T_{L}(\tau)$ for each $\tau \in \sigma^{-1}(M)$, $L^M$ gives a shelling of a subcomplex of $skel_d(\Lambda)$ with $h$-vector equal to the $F$-vector of $M$.

We will need to restrict our attention to a special class of multicomplexes. Define a partial order on our monomials as follows. For $\mu= x_1^{c_1}x_2^{c_2}\ldots x_k^{c_k}$ and $\mu'=x_1^{d_1}x_2^{d_2}\ldots x_k^{d_k}$ elements of $S(a_1, a_2, \ldots , a_k)$ with deg($\mu$) = deg($\mu'$), say $\mu < \mu'$ if for some $i$, $c_i < d_i$ and $c_j = d_j$ for all $j > i$ (reverse lexicographical order within degrees). For $F = (F_0, F_1, \ldots)$ the $F$-vector of some multicomplex $M \subseteq S(a_1, a_2, \ldots , a_k)$, let $S_{i,F_i}$ be the set of the first $F_i$ degree $i$ elements of $S(a_1, \ldots , a_k)$ in the reverse lex order, and set $I_M = \cup_{i \geq 0} S_{i,F_i}$. A result of Clements and Lindstr\"{o}m will allow us to replace $M$ with $I_M$:

\begin{theorem} \cite{MR0246781} Suppose $M$ is a multicomplex in $S(a_1, a_2, \ldots , a_k)$, where $a_1 \geq a_2 \geq \ldots \geq a_k$. Then $I_M$ is a multicomplex.
\end{theorem}

In particular, we may from now on assume that our multicomplex $M$ has the property that if deg($\mu$) = deg($\mu'$), $\mu < \mu'$ and $\mu' \in M$, then $\mu \in M$ (as $I_M$ clearly has this property and $F(I_M) = F(M)$). Thus, it will suffice to construct $L$ and $\sigma$ such that whenever $\gamma \in T_L(\tau_i)$, there exists $j < i$ and divisor $\mu$ of $\sigma(\tau_i)$ such that $\gamma \subseteq \tau_j$, $\text{deg}(\mu) = \text{deg}(\sigma(\tau_j))$, and $\sigma(\tau_j) \leq \mu$. Then if $\tau_i \in \sigma^{-1}(M)$, the properties of $M$ require that $\sigma(\tau_j) \in M$, so $\tau_j \in \sigma^{-1}(M)$, and then as $\gamma \subseteq \tau_j$, $T_{L^M}(\tau_i) = T_{L}(\tau_i)$.

\section{An Illustrative Example}
\label{ilex}
At this point it will be helpful to look at a small but non-trivial example. Let $d=4$ and $\Lambda = \Lambda(0;3,3) = \partial \Delta^2 * \partial \Delta^2$. The vertex set $V$ of $\Lambda$ decomposes into the vertex sets $P_1$ and $P_2$ of the two copies of $\partial \Delta^2$. The faces of $skel_4(\Lambda)$ are precisely the subsets of $V$ of size $4$ that do not contain either $P_1$ or $P_2$. Label the vertices of $\Lambda$ as shown:

\begin{figure}[h] 
\setlength{\unitlength}{0.254mm}
\begin{picture}(102,76)(45,-88)

        \allinethickness{0.254mm}\special{sh 0.3}\put(50,-45){\ellipse{4}{4}} 
        \allinethickness{0.254mm}\special{sh 0.3}\put(50,-60){\ellipse{4}{4}} 
        \allinethickness{0.254mm}\special{sh 0.3}\put(50,-75){\ellipse{4}{4}} 
        \allinethickness{0.254mm}\special{sh 0.3}\put(110,-45){\ellipse{4}{4}} 
        \allinethickness{0.254mm}\special{sh 0.3}\put(110,-60){\ellipse{4}{4}} 
        \allinethickness{0.254mm}\special{sh 0.3}\put(110,-75){\ellipse{4}{4}} 
        \put(45,-26){\textbf{\shortstack{$P_1$}}} 
        \put(110,-26){\textbf{\shortstack{$P_2$}}} 
        \put(60,-46){\shortstack{1}} 
        \put(120,-46){\shortstack{2}} 
        \put(60,-66){\shortstack{3}} 
        \put(120,-66){\shortstack{4}} 
        \put(60,-86){\shortstack{5}} 
        \put(120,-86){\shortstack{6}} 
\end{picture}
\caption{Vertex set of $\Lambda(0;3,3)$}
\label{exfig}
\end{figure}

We want to build a shelling of $skel_4(\Lambda)$ and a correspondence $\sigma$ between the facets of $skel_4(\Lambda)$ and the elements of $S(2,2)$ with the properties described at the end of the last section. Given our use of the reverse lexicographical order on the set of monomials, it is tempting to simply list the facets in reverse lex order $L_R$ (which will indeed give a shelling) and for $\tau$ the $i^{\text{th}}$ facet of $skel_4(\Lambda)$ having $|T_{L_R}(\tau)| = j$, let $\sigma(\tau)$ be the $i^{\text{th}}$ monomial in $S(2,2)$ of degree $j$. In fact such an approach will work in some simple cases. Here, however, it fails: 
\begin{center}
\begin{tabular}{c | c | c || c | c | c}
$\tau$ & $|T_{L_R}(\tau)|$ & $\sigma(\tau)$ & $\tau$ & $|T_{L_R}(\tau)|$ & $\sigma(\tau)$\\ \hline
1234 & 0 &1 & 1256 & 2 & $x_2^2$\\
1245 & 1 & $x_1$ &2356 & 3 & $x_1^2x_2$\\
2345 & 2 & $x_1^2$ & 1456 & 3 & $x_1x_2^2$\\
1236 & 1 & $x_2$ & 3456 & 4 & $x_1^2x_2^2$\\
1346 & 2 & $x_1x_2$ & & &  \\
\end{tabular}
\end{center}

In particular, consider the multicomplex $M = \{ 1, x_1, x_2, x_1^2, x_1x_2, x_1^2x_2 \}$. Note that $M = I_M$, but $L^M = (1234, 1245, 2345, 1236, 1346, 2356$). Then $T_{L^M}(2356) = \{ 235, 236 \} \neq T_L(2356)$, and letting $\Gamma = \cup_{\tau \in \sigma^{-1}(M)}\overline{\tau}$, $h(\Gamma) = (1,2,3,0, 0) \neq F(M)$. The problem is that $T_L(2356) = \{ 235, 236, 256 \}$, and these faces first appear in facets corresponding to $x_1^2, x_2$ and $x_2^2$. But $\sigma(2356) = x_1^2x_2$, the presence of which in $M$ does not imply that of $x_2^2$.

Let us examine the problem more closely. Notice that our ordering on the vertex set has resulted in each facet ending in $5$ corresponding to a monomial with greatest variable $x_1$, and any facet ending in $6$ corresponding to a monomial with greatest variable $x_2$. This leads us to define $\Lambda_i = \{ \gamma \in link_{\Lambda}(y_i) : \gamma \subseteq \{y_1, \ldots , y_{i-1} \} \}$ and $S_i = \{ \mu \in S : \text{supp}(\mu) \subseteq \{x_1, \ldots , x_i \} \text{ and } \mu x_i \in S\}$,  with the observation that any facet of $skel_d(\Lambda)$ is, for some $i$, of the form $\gamma \cup y_i$, where $\gamma \in skel_{d-1}(\Lambda_i)$, and any element of $S^d$ (aside from 1) is, for some $i$, of the form $\mu x_i$, where $\mu \in S^{d-1}_i$.

Consider $\Lambda_6$. This is isomorphic to $\Lambda(0;3, 2)$.
Note that our original ordering of facets gives a shelling of $skel_3(\Lambda_6)$ and correspondence $\sigma'$ to elements of $S(2,1)$, by taking $\sigma'(\tau) = \frac{\sigma(\tau \cup \{ 6 \})}{x_2}$.

\begin{center}
  \begin{tabular}{c | c | c}
  $\tau$ & $|T(\tau)|$ & $\sigma'(\tau)$\\ \hline
  123 & 0& 1\\
  134 & 1& $x_1$\\
  125 & 1& $x_2$\\
  235 & 2& $x_1^2$\\
  145 & 2& $x_1x_2$\\
  345 & 3& $x_1^2x_2$\\
  \end{tabular}
\end{center}

Here we see the same problem as before, occurring at 235. Na\"{i}vely we might note that here we no longer have the nice correspondence between last variable and last vertex we had in the larger ordering, but this deficiency is easily fixed by a simple reordering of the vertex set. In fact, consider the shelling and map obtained if we order our facets as if $ 4 >5$, while retaining our ordering on the monomials:

\begin{center}
  \begin{tabular}{c | c | c}
  $\tau$ & $|T(\tau)|$ & $\sigma'(\tau)$\\ \hline
  123 & 0& 1\\
  125 & 1& $x_1$\\
  235 & 2& $x_1^2$\\
  134 & 1& $x_2$\\
  154 & 2& $x_1x_2$\\
  354 & 3& $x_1^2x_2$\\
  \end{tabular}
\end{center}

It is simple to check that this correspondence has the property described at the end of the previous section, and we may furthermore, we can use this to fix our original attempt, by reordering the facets ending in 6 to match our new ordering on the facets of $skel_3(\Lambda_6)$:

\begin{center}
\begin{tabular}{c | c | c || c | c | c}
$\tau$ & $|T(\tau)|$ & $\sigma(\tau)$ & $\tau$ & $|T(\tau)|$ & $\sigma(\tau)$\\ \hline
1234 & 0 &1 & 2356 & 3 & $x_1^2x_2$\\
1245 & 1 & $x_1$ &1346 & 2 & $x_2^2$\\
2345 & 2 & $x_1^2$ & 1456 & 3 & $x_1x_2^2$\\
1236 & 1 & $x_2$ & 3456 & 4 & $x_1^2x_2^2$\\
1256 & 2 & $x_1x_2$ & & &  \\
\end{tabular}
\end{center}

The example suggests that we should build our shelling and map $\sigma$ inductively, at each step making sure the vertices are ordered so that the last $m$ vertices are from $P_1, P_2, \ldots P_m$, respectively. This is how we shall proceed. 

\section{Construction of the Shelling and Bijection}
\label{construct}
Let $\Lambda = \Lambda(l; p_1, \ldots , p_m)$ with $p_1 \geq p_2 \geq \ldots \geq p_m$. Let $V'$ be the vertex set of the $\Delta^{l-1}$ in the construction of $\Lambda$ and for $1 \leq i \leq m$ let $P_i$ be the vertex set of $\partial \Delta^{p_i - 1}$. (We will now allow $p_i =1$, in which case $P_i = \emptyset$, for the sake of an induction argument to come; similarly, we will allow $S(a_1, \ldots a_k)$ where $a_i = 0$, in which case we simply drop the variable $x_i$). Let $V = V' \cup (\cup_i P_i)$, and let $n = l + \sum |P_i|$. Let $S = S(\infty^{n-d-m},p_1-1, p_2-1, \ldots , p_m-1)$.

 As we will be changing the ordering on the vertices at different steps of our induction, we will require some additional notation. For $O$ denoting a total ordering $y_1 < y_2 < \ldots < y_n$ of $V$, let $\Lambda_{k,O}$ be $\Lambda_k$, as defined in the previous section, with respect to ordering $O$. (The ordering $x_1, x_2, \ldots x_{n-d}$ will remain fixed, so $S_k$ may remain as above.) 
 
 Recall that one characterization of a shelling $L = (\tau_1, \tau_2, \ldots , \tau_r)$ is that for each $i$ there exists a face $R(\tau_i)$ of $\tau_i$ such that $\overline{\tau_i} - (\cup_{j<i}\overline{\tau_j}) = \{ \gamma \subseteq \tau_i : R(\tau_i) \subseteq \gamma \}$. (Note in particular that $|T_L(\tau_i)| = |R_i|$.) Examining the two shellings of $skel_4(\Lambda(0;3,3))$ in our example in the last section, we see that both yield the same $R(\tau)$ for each facet $\tau$ of $skel_4(\Lambda)$. It will be helpful to determine the exact structure of the $R(\tau)$ in the shelling obtained by listing the facets  of $skel_d(\Lambda)$ in the reverse lexicographical order.
 
 Let $\tau$ be a face of $\Lambda$, and $O$ an ordering of $V$. Then let $full(\tau) = \{ i : |P_i \cap \tau| = |P_i|-1 \}$, and for $i \in full(\tau)$ let $miss(\tau,i)$ be the element of $P_i$ \emph{not} in $\tau$ (the notation is meant to suggest that $full(\tau)$ collects the indices of the sets $P_i$ such that $\tau \cap P_i$ is `full' in the sense that no further elements of $P_i$ could be added without leaving $\Lambda$, and $miss(\tau, i)$ is the element of $P_i$ missing from $\tau$). Let $s_O(\tau)$ be the first element of $V - \cup_{i \in full(\tau)} \{ miss(\tau,i) \}$ not appearing in $\tau$ (with respect to order $O$) if such an element exists, otherwise set $s_O(\tau) = \infty$. Let $\tau_{>s_O} = \{y \in \tau : y > s_O \}$,  and $U_O(\tau) = \{ y :  y \in P_i \text{ and } y>miss(\tau,i) \text{ for some } i\in full(\tau) \}$. Finally, let $R_O(\tau) = \tau_{>(s_O(\tau))} \cup U_O(\tau)$.
 
 \begin{example} \label{uex} Let $\Lambda = \Lambda(1;5,4,3)$, with vertex ordering $O$ as shown:
 \begin{figure}[h]
\setlength{\unitlength}{0.254mm}
\begin{picture}(151,91)(20,-98)

        \allinethickness{0.254mm}\special{sh 0.3}\put(55,-30){\ellipse{4}{4}} 
        \allinethickness{0.254mm}\special{sh 0.3}\put(80,-30){\ellipse{4}{4}} 
        \allinethickness{0.254mm}\special{sh 0.3}\put(80,-45){\ellipse{4}{4}} 
        \allinethickness{0.254mm}\special{sh 0.3}\put(80,-60){\ellipse{4}{4}} 
        \put(20,-21){\textbf{\shortstack{$V'$}}} 
        \put(50,-21){\textbf{\shortstack{$P_1$}}} 
        \put(75,-21){\textbf{\shortstack{$P_2$}}} 
        \allinethickness{0.254mm}\special{sh 0.3}\put(80,-75){\ellipse{4}{4}} 
        \put(100,-21){\textbf{\shortstack{$P_3$}}} 
        \allinethickness{0.254mm}\special{sh 0.3}\put(25,-30){\ellipse{4}{4}} 
        \allinethickness{0.254mm}\special{sh 0.3}\put(55,-45){\ellipse{4}{4}} 
        \allinethickness{0.254mm}\special{sh 0.3}\put(55,-60){\ellipse{4}{4}} 
        \allinethickness{0.254mm}\special{sh 0.3}\put(55,-75){\ellipse{4}{4}} 
        \allinethickness{0.254mm}\special{sh 0.3}\put(105,-30){\ellipse{4}{4}} 
        \allinethickness{0.254mm}\special{sh 0.3}\put(105,-45){\ellipse{4}{4}} 
        \allinethickness{0.254mm}\special{sh 0.3}\put(105,-60){\ellipse{4}{4}} 
        \put(110,-66){\textbf{\shortstack{$y_{13}$}}} 
        \put(110,-51){\textbf{\shortstack{$y_7$}}} 
        \put(110,-36){\textbf{\shortstack{$y_4$}}} 
        \put(85,-66){\textbf{\shortstack{$y_9$}}} 
        \put(85,-51){\textbf{\shortstack{$y_6$}}} 
        \put(85,-36){\textbf{\shortstack{$y_3$}}} 
        \put(85,-81){\textbf{\shortstack{$y_{12}$}}} 
        \put(60,-96){\textbf{\shortstack{$y_{11}$}}} 
        \allinethickness{0.254mm}\special{sh 0.3}\put(55,-90){\ellipse{4}{4}} 
        \put(60,-81){\textbf{\shortstack{$y_{10}$}}} 
        \put(60,-66){\textbf{\shortstack{$y_8$}}} 
        \put(60,-51){\textbf{\shortstack{$y_5$}}} 
        \put(60,-36){\textbf{\shortstack{$y_2$}}} 
        \put(30,-36){\textbf{\shortstack{$y_1$}}} 
\end{picture}
\caption{Vertex set of $\Lambda(1;5,4,3)$ with ordering $O$}
 \end{figure}
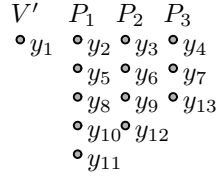

 Consider the face $\tau = \{y_1,y_2,y_4,y_5,y_6,y_9, y_{11}, y_{12} \}$. Then $full(\tau) = \{ 2 \}$, $miss(\tau, 2) = y_3$, $U_O(\tau) = \{ y_6, y_9, y_{12} \}$, $s_O(\tau) = y_7$, and $\tau_{>s_O(\tau)} = \{ y_9, y_{11}, y_{12} \}$.  So $R_O(\tau) = \{ y_6, y_9, y_{11}, y_{12} \}$.
 
 \begin{figure}[h]
 \begin{minipage}{.3\linewidth}
 \centering
\setlength{\unitlength}{0.254mm}
\begin{picture}(161,93)(20,-100)

        \allinethickness{0.254mm}\special{sh 0.3}\put(55,-30){\ellipse{4}{4}} 
        \allinethickness{0.254mm}\special{sh 0.3}\put(85,-30){\ellipse{4}{4}} 
        \allinethickness{0.254mm}\special{sh 0.3}\put(85,-45){\ellipse{4}{4}} 
        \allinethickness{0.254mm}\special{sh 0.3}\put(85,-60){\ellipse{4}{4}} 
        \put(20,-21){\textbf{\shortstack{$V'$}}} 
        \put(50,-21){\textbf{\shortstack{$P_1$}}} 
        \put(80,-21){\textbf{\shortstack{$P_2$}}} 
        \allinethickness{0.254mm}\special{sh 0.3}\put(85,-75){\ellipse{4}{4}} 
        \put(110,-21){\textbf{\shortstack{$P_3$}}} 
        \allinethickness{0.254mm}\special{sh 0.3}\put(25,-30){\ellipse{4}{4}} 
        \allinethickness{0.254mm}\special{sh 0.3}\put(55,-45){\ellipse{4}{4}} 
        \allinethickness{0.254mm}\special{sh 0.3}\put(55,-60){\ellipse{4}{4}} 
        \allinethickness{0.254mm}\special{sh 0.3}\put(55,-75){\ellipse{4}{4}} 
        \allinethickness{0.254mm}\special{sh 0.3}\put(115,-30){\ellipse{4}{4}} 
        \allinethickness{0.254mm}\special{sh 0.3}\put(115,-45){\ellipse{4}{4}} 
        \allinethickness{0.254mm}\special{sh 0.3}\put(115,-60){\ellipse{4}{4}} 
        \put(120,-66){\textbf{\shortstack{$y_{13}$}}} 
        \put(120,-51){\textbf{\shortstack{$y_7$}}} 
        \put(120,-36){\textbf{\shortstack{$y_4$}}} 
        \put(90,-66){\textbf{\shortstack{$y_9$}}} 
        \put(90,-51){\textbf{\shortstack{$y_6$}}} 
        \put(90,-36){\textbf{\shortstack{$y_3$}}} 
        \put(90,-81){\textbf{\shortstack{$y_{12}$}}} 
        \put(60,-96){\textbf{\shortstack{$y_{11}$}}} 
        \allinethickness{0.254mm}\special{sh 0.3}\put(55,-90){\ellipse{4}{4}} 
        \put(60,-81){\textbf{\shortstack{$y_{10}$}}} 
        \put(60,-66){\textbf{\shortstack{$y_8$}}} 
        \put(60,-51){\textbf{\shortstack{$y_5$}}} 
        \put(60,-36){\textbf{\shortstack{$y_2$}}} 
        \put(30,-36){\textbf{\shortstack{$y_1$}}} 
        \allinethickness{0.254mm}\path(110,-40)(110,-85)(80,-85)(80,-55)(50,-55)(50,-40)(20,-40)(20,-25)(80,-25)(80,-40)(110,-40) 
        \allinethickness{0.254mm}\path(50,-85)(80,-85)(80,-100)(50,-100)(50,-85) 
        \allinethickness{0.254mm}\path(110,-25)(110,-35)(110,-40)(140,-40)(140,-25)(110,-25) 
\end{picture}
\caption{$\tau$}
 \end{minipage}%
 \begin{minipage}{.3\linewidth}
\setlength{\unitlength}{0.254mm}
\begin{picture}(161,91)(20,-98)

        \allinethickness{0.254mm}\special{sh 0.3}\put(55,-30){\ellipse{4}{4}} 
        \allinethickness{0.254mm}\special{sh 0.3}\put(85,-30){\ellipse{4}{4}} 
        \allinethickness{0.254mm}\special{sh 0.3}\put(85,-45){\ellipse{4}{4}} 
        \allinethickness{0.254mm}\special{sh 0.3}\put(85,-60){\ellipse{4}{4}} 
        \put(20,-21){\textbf{\shortstack{$V'$}}} 
        \put(50,-21){\textbf{\shortstack{$P_1$}}} 
        \put(80,-21){\textbf{\shortstack{$P_2$}}} 
        \allinethickness{0.254mm}\special{sh 0.3}\put(85,-75){\ellipse{4}{4}} 
        \put(110,-21){\textbf{\shortstack{$P_3$}}} 
        \allinethickness{0.254mm}\special{sh 0.3}\put(25,-30){\ellipse{4}{4}} 
        \allinethickness{0.254mm}\special{sh 0.3}\put(55,-45){\ellipse{4}{4}} 
        \allinethickness{0.254mm}\special{sh 0.3}\put(55,-60){\ellipse{4}{4}} 
        \allinethickness{0.254mm}\special{sh 0.3}\put(55,-75){\ellipse{4}{4}} 
        \allinethickness{0.254mm}\special{sh 0.3}\put(115,-30){\ellipse{4}{4}} 
        \allinethickness{0.254mm}\special{sh 0.3}\put(115,-45){\ellipse{4}{4}} 
        \allinethickness{0.254mm}\special{sh 0.3}\put(115,-60){\ellipse{4}{4}} 
        \put(120,-66){\textbf{\shortstack{$y_{13}$}}} 
        \put(120,-51){\textbf{\shortstack{$y_7$}}} 
        \put(120,-36){\textbf{\shortstack{$y_4$}}} 
        \put(90,-66){\textbf{\shortstack{$y_9$}}} 
        \put(90,-51){\textbf{\shortstack{$y_6$}}} 
        \put(90,-36){\textbf{\shortstack{$y_3$}}} 
        \put(90,-81){\textbf{\shortstack{$y_{12}$}}} 
        \put(60,-96){\textbf{\shortstack{$y_{11}$}}} 
        \allinethickness{0.254mm}\special{sh 0.3}\put(55,-90){\ellipse{4}{4}} 
        \put(60,-81){\textbf{\shortstack{$y_{10}$}}} 
        \put(60,-66){\textbf{\shortstack{$y_8$}}} 
        \put(60,-51){\textbf{\shortstack{$y_5$}}} 
        \put(60,-36){\textbf{\shortstack{$y_2$}}} 
        \put(30,-36){\textbf{\shortstack{$y_1$}}} 
        \allinethickness{0.254mm}\path(80,-40)(110,-40)(110,-85)(80,-85)(80,-40) 
\end{picture}
\centering
\caption{$U_O(\tau)$}
 \end{minipage}%
  \begin{minipage}{.3\linewidth}
\setlength{\unitlength}{0.254mm}
\begin{picture}(161,93)(20,-100)

        \allinethickness{0.254mm}\special{sh 0.3}\put(55,-30){\ellipse{4}{4}} 
        \allinethickness{0.254mm}\special{sh 0.3}\put(85,-30){\ellipse{4}{4}} 
        \allinethickness{0.254mm}\special{sh 0.3}\put(85,-45){\ellipse{4}{4}} 
        \allinethickness{0.254mm}\special{sh 0.3}\put(85,-60){\ellipse{4}{4}} 
        \put(20,-21){\textbf{\shortstack{$V'$}}} 
        \put(50,-21){\textbf{\shortstack{$P_1$}}} 
        \put(80,-21){\textbf{\shortstack{$P_2$}}} 
        \allinethickness{0.254mm}\special{sh 0.3}\put(85,-75){\ellipse{4}{4}} 
        \put(110,-21){\textbf{\shortstack{$P_3$}}} 
        \allinethickness{0.254mm}\special{sh 0.3}\put(25,-30){\ellipse{4}{4}} 
        \allinethickness{0.254mm}\special{sh 0.3}\put(55,-45){\ellipse{4}{4}} 
        \allinethickness{0.254mm}\special{sh 0.3}\put(55,-60){\ellipse{4}{4}} 
        \allinethickness{0.254mm}\special{sh 0.3}\put(55,-75){\ellipse{4}{4}} 
        \allinethickness{0.254mm}\special{sh 0.3}\put(115,-30){\ellipse{4}{4}} 
        \allinethickness{0.254mm}\special{sh 0.3}\put(115,-45){\ellipse{4}{4}} 
        \allinethickness{0.254mm}\special{sh 0.3}\put(115,-60){\ellipse{4}{4}} 
        \put(120,-66){\textbf{\shortstack{$y_{13}$}}} 
        \put(120,-51){\textbf{\shortstack{$y_7$}}} 
        \put(120,-36){\textbf{\shortstack{$y_4$}}} 
        \put(90,-66){\textbf{\shortstack{$y_9$}}} 
        \put(90,-51){\textbf{\shortstack{$y_6$}}} 
        \put(90,-36){\textbf{\shortstack{$y_3$}}} 
        \put(90,-81){\textbf{\shortstack{$y_{12}$}}} 
        \put(60,-96){\textbf{\shortstack{$y_{11}$}}} 
        \allinethickness{0.254mm}\special{sh 0.3}\put(55,-90){\ellipse{4}{4}} 
        \put(60,-81){\textbf{\shortstack{$y_{10}$}}} 
        \put(60,-66){\textbf{\shortstack{$y_8$}}} 
        \put(60,-51){\textbf{\shortstack{$y_5$}}} 
        \put(60,-36){\textbf{\shortstack{$y_2$}}} 
        \put(30,-36){\textbf{\shortstack{$y_1$}}} 
        \allinethickness{0.254mm}\path(50,-85)(80,-85)(80,-100)(50,-100)(50,-85) 
        \allinethickness{0.254mm}\path(80,-55)(110,-55)(110,-85)(80,-85)(80,-55) 
\end{picture}

\centering
\caption{$\tau_{>s_O(\tau)}$}
 \end{minipage}%
  \end{figure}
 \end{example}
 
 Now, if $\tau$ is a facet of $skel_d(\Lambda)$ and $\gamma$ is a facet of $\overline{\tau}$,  $\gamma = \tau - y_i$ for some $y_i \in \tau$. Then $\gamma$ appears as a face of a facet occurring before $\tau$ in the reverse lex order (as determined by the ordering $O$ on the vertices) if and only if $\gamma \cup y_j$ is a facet of $skel_d(\Lambda)$ for some $j < i$. It may easily be checked that $\tau$ is the reverse lexicographically first facet of $skel_d(\Lambda)$ containing $R_O(\tau)$, so if $R_O(\tau) \subseteq \gamma$, $\gamma$ occurs in no earlier facet. On the other hand, if there is $y_r$ in $R_O(\tau)$ such that $y_r \notin \gamma$, either $y_r \in U_O(\tau)$, in which case $\gamma \cup miss(\tau, k)$  (where $y_r \in P_k$) is a reverse lexicographically earlier facet of $skel_d(\Lambda)$ containing $\gamma$, or $y_r > s_O(\tau)$, in which case $\gamma \cup s_O(\tau)$ is an earlier facet of $skel_d(\Lambda)$ containing $\gamma$. Thus, if $L =( \tau_1, \tau_2, \ldots, \tau_r)$ is the reverse lex order on the facets of $skel_d(\Lambda)$, $\overline{\tau_i} - (\cup_{j<i}\overline{\tau_j}) = \{ \gamma \subseteq \tau : R_O(\tau_i) \subseteq \gamma \}$. Our inductively built shelling will share this structure. 
 
 We are now ready to prove our central theorem.

\begin{theorem} \label{big} Let $\Lambda$ and $S$ be as above, and let $O$ be an ordering $y_1 < y_2 < \ldots < y_n$ of $V$ such that for $1 \leq i \leq m$, $y_{n-m+i} \in P_i$. Let $1 \leq d \leq n-m$. Then there exists a shelling $L = (\tau_1, \tau_2, \ldots ,\tau_t)$ of $skel_d(\Lambda)$ and bijection $\sigma$ from the set of facets of $skel_d(\Lambda)$ to $S^d$ such that:
\begin{enumerate}
\item \label{R} $\overline{\tau_i} - (\cup_{j<i}\overline{\tau_j}) = \{ \gamma \subseteq \tau : R_O(\tau) \subseteq \gamma \}$.
\item \label{deg} $deg(\sigma(\tau_i)) = |T_L(\tau_i)|$.
\item \label{ord} If $\gamma \in T_L(\tau_i)$, then there exists $j <i$ and $\mu | \sigma(\tau_i)$ such that $\gamma \subseteq \tau_j$, deg($\mu$) = deg($\sigma(\tau_j)$), and $\sigma(\tau_j) \leq \mu$.

\end{enumerate}
\end{theorem}
Again, condition (\ref{R}) is sufficient to show that $L$ is a shelling. Theorem \ref{main} follows from (2) and (3), the proof of the latter requiring our precise definition of $R_O$.
\begin{proof}  (of Theorem \ref{big})

We will proceed by induction on $d$. If $d =1$, let $L = (y_1, y_2, \ldots y_n$), $\sigma(y_1) = 1$, and $\sigma(y_i) = x_{i-1}$ for $1 < i \leq n$. Properties (1)-(3) immediately follow.

Now, suppose $1 <d \leq n-m$. Set $\tau^0_1 = \{ y_1, y_2, \ldots , y_d \}$. By the properties of our order on $V$, $\tau^0_1$ does not contain any $P_i$, and hence is a facet of $skel_d(\Lambda)$. Set $\sigma(\tau^0_1) = 1$. 

Any other facet $\tau$ of $skel_d(\Lambda)$ has the form $\tau = \tau' \cup y_{d+k}$ where $\tau' \in skel_{d-1}(\Lambda_{d+k, O})$ for some $k>0$. Similarly, any element of $S^d$ aside from $1$ is of the form $\mu x_k$ where $\mu \in S^{d-1}_k$ for some $k \geq 1$.

Suppose $d+k \leq n-m$. Then $skel_{d-1}(\Lambda_{d+k, O})$ is simply the ($d-2$)-skeleton of the simplex on the first $d+k-1$ vertices in $V$. Then the ordering $O_k$ on these vertices inherited from the original order on $V$ satisfies the condition of our theorem, so by induction there exists a shelling of $skel_d(\Lambda_{d+k, O_k})$ and map $\sigma_k$ from its set of facets to $S^{d-1}(\infty^{d+k-1 - (d-1)}) = S^{d-1}(\infty^k) = S^{d-1}_k$ satisfying (1)-(3). Call this shelling $L_k = (G^k_1, G^k_2, \ldots, G^k_{r_k})$.

On the other hand, if $d+k = n-m +i$ for some $1 \leq i \leq m$, then $skel_{d-1}(\Lambda_{d+k,O}) = skel_{d-1}(\Lambda(l + \sum_{j >i} (p_i-1); p_1, p_2, \ldots, p_{i-1}, (p_i - 1)))$. In this case, the restriction of the order on $V$ does not quite meet the conditions of the theorem. Let $y^k$ be the largest element of $P_i - y_{d+k}$ with respect to $O$, and define a new order $O_k$ by $y_1 <_k y_2 <_k \ldots <_k \widehat{y^k} <_k \ldots <_k y_{d+k-1} <_k y^k$ ( i.e., take the original order but set $y_i <_k y^k$ for all $y_i \neq y^k$, as in the example in the previous section).  This new order satisfies the conditions of our theorem, and so by induction we have a shelling $L_k = (G^k_1, G^k_2, \ldots, G^k_{r_k})$ of $skel_{d-1}(\Lambda_{d+k, O_k})$ and map $\sigma_k$ from the set of its facets to $S^{d-1}(\infty^{d+k-1-(d-1)-i}, p_1-1, \ldots , p_{i-1}-1, p_i-2) = S^{d-1}(\infty^{n-m-d}, p_1-1, \ldots , p_{i-1}-1, p_i-2)= S^{d-1}_{k}$ satisfying (1) - (3).

For $1 \leq k \leq n-d$, and $1 \leq i \leq r_k$, set $\tau^k_i = G^k_i \cup y_{d+k}$ and $\sigma(\tau^k_i) = \sigma(G^k_i)x_k$. Let $L = (\tau^0_1, \tau^1_1, \ldots , \tau^1_{r_1}, \tau^2_1, \ldots, \tau^2_{r_2}, \ldots , t^{n-d}_{r_{n-d}})$. Our claim then is that $L$ and $\sigma$ satisfy (1)-(3).

(1) The $k=0$ case is immediate. Now suppose $k >0$. Set $R_i^k = R_{O_k}(G_i^k) \cup y_{d+k}$. We will first show that $\overline{\tau^k_i} - (\cup_{k' < k}\overline{\tau^{k'}_j} \cup (\cup_{j<i} \overline{\tau^k_j})) = \{ \gamma \subseteq \tau_i^k : R_i^k \subseteq \gamma \}$, and then that $R_i^k = R_O(\tau_i^k$).

\begin{example} Let $\Lambda$, $O$, and $\tau$ be as in Example \ref{uex}. Then $\tau = G \cup y_{12}$, where $G \in skel_7(\Lambda_{12}) = skel_7(\Lambda(3;5,3))$. In the new ordering $O_{12}$, $y_9$ becomes the last vertex, so labeling the $i^{th}$ vertex in this ordering $y'_i$, we have $y'_i = y_i$ for $i < 9$, $y'_9 = y_{10}$, $y'_{10} = y_{11}$, and $y'_{11} = y_9$.
Observe that $R_O(\tau) = R_{O_{12}}(G) \cup y_{12}$.
\begin{figure}[h]
\begin{minipage}{.5\linewidth}
\centering
\setlength{\unitlength}{0.254mm}
\begin{picture}(131,93)(20,-100)

        \allinethickness{0.254mm}\special{sh 0.3}\put(55,-30){\ellipse{4}{4}} 
        \allinethickness{0.254mm}\special{sh 0.3}\put(85,-30){\ellipse{4}{4}} 
        \allinethickness{0.254mm}\special{sh 0.3}\put(85,-45){\ellipse{4}{4}} 
        \allinethickness{0.254mm}\special{sh 0.3}\put(85,-60){\ellipse{4}{4}} 
        \put(20,-21){\textbf{\shortstack{$V'$}}} 
        \put(50,-21){\textbf{\shortstack{$P_1$}}} 
        \put(80,-21){\textbf{\shortstack{$P_2$}}} 
        \allinethickness{0.254mm}\special{sh 0.3}\put(25,-30){\ellipse{4}{4}} 
        \allinethickness{0.254mm}\special{sh 0.3}\put(55,-45){\ellipse{4}{4}} 
        \allinethickness{0.254mm}\special{sh 0.3}\put(55,-60){\ellipse{4}{4}} 
        \allinethickness{0.254mm}\special{sh 0.3}\put(55,-75){\ellipse{4}{4}} 
        \allinethickness{0.254mm}\special{sh 0.3}\put(25,-45){\ellipse{4}{4}} 
        \allinethickness{0.254mm}\special{sh 0.3}\put(25,-60){\ellipse{4}{4}} 
        \put(30,-66){\textbf{\shortstack{$y'_7$}}} 
        \put(30,-51){\textbf{\shortstack{$y'_4$}}} 
        \put(90,-66){\textbf{\shortstack{$y'_{11}$}}} 
        \put(90,-51){\textbf{\shortstack{$y'_6$}}} 
        \put(90,-36){\textbf{\shortstack{$y'_3$}}} 
        \put(60,-96){\textbf{\shortstack{$y'_{10}$}}} 
        \allinethickness{0.254mm}\special{sh 0.3}\put(55,-90){\ellipse{4}{4}} 
        \put(60,-81){\textbf{\shortstack{$y'_9$}}} 
        \put(60,-66){\textbf{\shortstack{$y'_8$}}} 
        \put(60,-51){\textbf{\shortstack{$y'_5$}}} 
        \put(60,-36){\textbf{\shortstack{$y'_2$}}} 
        \put(30,-36){\textbf{\shortstack{$y'_1$}}} 
        \allinethickness{0.254mm}\path(20,-25)(75,-25)(75,-40)(110,-40)(110,-70)(80,-70)(80,-55)(20,-55)(20,-25) 
        \allinethickness{0.254mm}\path(50,-85)(50,-100)(80,-100)(80,-85)(50,-85) 
\end{picture}
\caption{$G$}
\end{minipage}%
\begin{minipage}{.5\linewidth}
\centering
\setlength{\unitlength}{0.254mm}
\begin{picture}(131,93)(20,-100)

        \allinethickness{0.254mm}\special{sh 0.3}\put(55,-30){\ellipse{4}{4}} 
        \allinethickness{0.254mm}\special{sh 0.3}\put(85,-30){\ellipse{4}{4}} 
        \allinethickness{0.254mm}\special{sh 0.3}\put(85,-45){\ellipse{4}{4}} 
        \allinethickness{0.254mm}\special{sh 0.3}\put(85,-60){\ellipse{4}{4}} 
        \put(20,-21){\textbf{\shortstack{$V'$}}} 
        \put(50,-21){\textbf{\shortstack{$P_1$}}} 
        \put(80,-21){\textbf{\shortstack{$P_2$}}} 
        \allinethickness{0.254mm}\special{sh 0.3}\put(25,-30){\ellipse{4}{4}} 
        \allinethickness{0.254mm}\special{sh 0.3}\put(55,-45){\ellipse{4}{4}} 
        \allinethickness{0.254mm}\special{sh 0.3}\put(55,-60){\ellipse{4}{4}} 
        \allinethickness{0.254mm}\special{sh 0.3}\put(55,-75){\ellipse{4}{4}} 
        \allinethickness{0.254mm}\special{sh 0.3}\put(25,-45){\ellipse{4}{4}} 
        \allinethickness{0.254mm}\special{sh 0.3}\put(25,-60){\ellipse{4}{4}} 
        \put(30,-66){\textbf{\shortstack{$y'_7$}}} 
        \put(30,-51){\textbf{\shortstack{$y'_4$}}} 
        \put(90,-66){\textbf{\shortstack{$y'_{11}$}}} 
        \put(90,-51){\textbf{\shortstack{$y'_6$}}} 
        \put(90,-36){\textbf{\shortstack{$y'_3$}}} 
        \put(60,-96){\textbf{\shortstack{$y'_{10}$}}} 
        \allinethickness{0.254mm}\special{sh 0.3}\put(55,-90){\ellipse{4}{4}} 
        \put(60,-81){\textbf{\shortstack{$y'_9$}}} 
        \put(60,-66){\textbf{\shortstack{$y'_8$}}} 
        \put(60,-51){\textbf{\shortstack{$y'_5$}}} 
        \put(60,-36){\textbf{\shortstack{$y'_2$}}} 
        \put(30,-36){\textbf{\shortstack{$y'_1$}}} 
        \allinethickness{0.254mm}\path(50,-85)(50,-100)(80,-100)(80,-85)(50,-85) 
        \allinethickness{0.254mm}\path(80,-40)(110,-40)(110,-70)(80,-70)(80,-40) 
\end{picture}
\caption{$R_{O_{12}}(G)$}
\end{minipage}
\end{figure}

\end{example}

Returning now to the proof, suppose $\gamma \subseteq \tau_i^k$ and $R_i^k \subseteq \gamma$. Then $y_{d+k} \in \gamma$, so $\gamma$ cannot be in any $\tau_i^{k'}$ for $k' < k$. On the other hand, as $\gamma - y_{d+k}$ contains $R_{O_k}(G_i^k)$, there is no $j<i$ such that $G_j^k$ contains $\gamma - y_{d+k}$. Hence $\gamma$ is not in any $\tau_j^k$ for $j<i$, so $\gamma$ can occur in no facet appearing before $\tau_i^k$.

Now suppose $R_i^k$ is not contained in $\gamma$. Then there is at least one element of $R_i^k$ not in $\gamma$. If some such element is in $R_{O_k}(G_i^k)$, then $\gamma - y_{d+k}$ is a face of $G_i^k$ not containing $R_{O_k}(G_i^k)$, so there is $j<i$ such that $\gamma - y_{d+k} \subseteq G_j^k$, and then $\gamma \subseteq \tau_j^k$. Otherwise, $\gamma = \tau_i^k - y_{d+k} = G_i^k$. Now, there is clearly some $r < d+k$ such that $y_r \notin G_i^k$. Suppose, in order to obtain a contradiction, that for each such $r$, $G_i^k \cup y_r$ is not a facet of $skel_d(\Lambda)$, i.e., $G_i^k \cup y_r$ contains some $P_s$. Then $d+k = n-m+j$ for some $j> 1$ (otherwise, $G_i^k$ cannot contain any element of the form $y_{n-m+s}$ for $1 \leq s \leq m$, and adding any vertex before $y_{d+k}$ cannot complete $P_s$). But then there are at least $|P_i| -1$ elements of each $P_j$ occurring before $y_{d+k}$ in our ordering, so $G_i^k$ must contain at least $|P_j| -1$ elements of each $P_j$, in addition to all of $V'$. Hence $d = |G_i^k| + 1 \geq l + \sum (|P_j|-1) + 1 \geq n-m+1$, a contradiction. Hence, there is some $r < d+k$ such that $G_i^k \cup y_r$ is a facet of $\Lambda^d$. This facet occurs before $\tau_i^k$ and contains $\gamma$.

It remains to show that $R_i^k = R_O(\tau_i^k)$ for $k>0$. We first confirm that $y_{d+k} \in R_O(\tau_i^k)$. If  $y_{d+k} > s_O(\tau_i^k)$, the $y_{d+k}$ is in $R_O(\tau_i^k)$, so suppose $y_{d+k} < s_O(\tau_i^k)$. Then as $y_{d+k}$ is the greatest element of $\tau_i^k$, $\tau_i^k$ must consist of every element of $V - \cup_{j \in full(\tau_j^k)} \{ miss(\tau_i^k,j) \}$ less than $y_{d+k}$. Suppose $d+k \leq n-m$. Then $\tau_i^k$ cannot contain the largest element of any $P_i$, so in particular  $miss(\tau_i^k,j) > y_{d+k}$ for any $j \in full(\tau_i^k)$. Thus $\tau_i^k$ is just the first $d$ elements of $V$, i.e. $\tau_1^0$. But $k >0$, so we must have $d+k = n-m+j$ for some $j$, and in particular $y_{d+k}$ is the largest element of $P_j$. But then $j \in full(\tau_i^k)$ and $y_{d+k} >  miss(\tau_i^k,j)$, so $y_{d+k} \in U_O(\tau_i^k) \subseteq R_O(\tau_i^k)$.


Now suppose $d+k \leq n-m$. In this case our orderings $O$ and $O_k$ are the same, so $s_{O_k}(G_i^k) = s_O(\tau_i^k)$. Furthermore, as $\tau_i^k$ cannot contain the largest element of any $P_j$,  $U_O(\tau_i^k) = \emptyset$. Thus $R_O(\tau_i^k) = R_{O_k}(G_i^k) \cup y_{d+k} = R_i^k$.

On the other hand, suppose $d+k = n-m+j$. Observe that the vertices corresponding to the indices in $full(G_i^k)$ are the same as those corresponding to the indices in $full(\tau_i^k)$, and as $y^r$ is the largest element of $P_j - y_{d+k}$ (with respect to both orders), with $O_k$ matching $O$ on all the other vertices, $U_O(\tau_i^k)$ and $U_{O_k}(G_i^k)$ agree, except for the possible presence of $y_{d+k}$ in the former. But we have already seen that $y_{d+k}$ must be in both $R_i^k$ and $R_O(\tau_i^k)$.

Suppose $s_{O_k}(G_i^k) < y^k$.  Then $s_O(\tau_i^k) = s_{O_k}(G_i^k)$. Furthermore, for $y \in G_i^k$, $y > s_O(\tau_i^k)$ if and only if $y >_k s_O(\tau_i^k)$. Thus $R_i^k = R_O(\tau_i^k)$.

On the other hand, suppose $s_{O_k}(G_i^k) \geq y^k$. As $y^k$ is the greatest element of $P_j - y_{d+k}$, any other element $y$ of $P_j - y_{d+k}$ is less than $s_{O_k}(G_i^k)$ in both orders. Thus, either all of these elements are in $G_i^k$ or exactly one is missing and every other element of $P_j - y_{d+k}$ is in $G_i^k$.

First, suppose every element of $P_j - y_{d+k}$ less than $y^k$ is in $G_i^k$. Then $y^k$ cannot be in $G_i^k$, and in particular $y^k = miss(G_i^k, j) = miss(\tau_i^k, j)$. Thus $s_O(\tau_i^k) = s_{O_k}(G_i^k)$ (as the changing of the position of $y^k$ in the order will have no effect on $s$). Furthermore, $s_O(\tau_i^k) \neq y^r$, and $y^r \notin \tau_i^k$. Then for $y \in G_i^k$, $y>_k s_{O_k}(G_i^k)$ if and only if $y > s_{O}(\tau_i^k)$. Thus $R_i^k = R_O(\tau_i^k)$.

Now suppose instead that $G_i^k$ contains every element of $P_j$ except some $y<y^r$. Then in particular $y^r \in G_i^k$, so $s_O(\tau_i^k) = s_{O_k}(G_i^k) \neq y^r$. Thus for $y \neq y^r$, $y>_k s_{O_k}(G_i^k)$ if and only if $y > s_{O}(\tau_i^k)$. It now only remains to check the membership of $y^r$ in $R_i^k$ and $R_O(\tau_i^k)$. But $y^r \in U_O(\tau_i^k) = U_{O_k}(G_i^k)$, and is thus in both $R_i^k$ and $R_O(\tau_i^k)$. Hence (1) is proved.
 
(2) Note that $|T_L(\tau_i^k)| = |R_O(\tau_i^k)| = |R_{O_k}(G_i^k)| + 1$. Then as $|R_{O_k}(G_i^k)| = \text{deg}(\sigma_k(G_i^k))$, (2) follows from the definition of $\sigma$.

(3) The $k=0$ case is trivial. Suppose $k>0$ and $\gamma \in T_L(\tau_i^k)$. Then $\gamma$ is obtained from $\tau_i^k$ by removing some element of $R_O(\tau_i^k)$. Suppose that element is not $y_{d+k}$. Then $\gamma - y_{d+k} \in T_{L_k}(G_i^k)$. Thus, there exists $j<i$ and divisor $\mu$ of $\sigma_k(G_i^k)$ such that $\gamma - y_{d+k} \in G_j^k$, deg$(\mu) = \text{deg}(\sigma_k(G_j^k))$, and $\sigma_k(G_j^k) \leq \mu$. Then $\gamma \in \tau^k_j$, $\mu x_k$ is a divisor of $\sigma(\tau_i^k)$, $\text{deg}(\mu x_k) = \text{deg}(\sigma(\tau_i^j)$, and $\sigma(\tau_i^j) \leq \mu x_k$.

On the other hand, suppose $\gamma = \tau_i^k - y_{d+k}$. We claim there exists a facet $\tau_t^r$ for some $r < k$ such that deg$(\sigma(\tau_t^r)) =$ deg$(\sigma(\tau_i^k))$.

First consider the case $y_{d+k} \in U_O(\tau_i^k)$, where $y_{d+k} \in P_j$. Then let $\tau' = \gamma \cup miss(\tau_i^k,j)$ . Note that $s_O(\tau') = s_O(\tau_i^k)$ and $full(\tau_i^k) = full(\tau')$.
Suppose $ y \in R_O(\tau')$ and $y\neq miss(\tau_i^k,j)$. Then $y \in \tau_i^k$, and if $y > s_O(\tau')$, $y > s_O(\tau_i^k)$, so $y \in R_O(\tau_i^k)$. On the other hand, if $y \in U_O(\tau')$, then $y \in U_O(\tau_i^k)$, as $miss(\tau',q) \geq miss(\tau_i^k, q)$ for all $q \in full(\tau') = full(\tau_i^k)$. Thus $y \in R_O(\tau_i^k)$. In particular, note that every element of $R_O(\tau') - miss(\tau_i^k,j)$ is in $R_O(\tau_i^k)$, and as $y_{d+k}$ is in $R_O(\tau_i^k)$ but not $\tau'$, $R_O(\tau') - miss(\tau_i^k,j) \subseteq R_O(\tau_i^k) - y_{d+k}$. Hence $|R_O(\tau')| \leq |R_O(\tau_i^k)|$.

\begin{example} Again take $\Lambda$, $O$, and $\tau$ as in Example \ref{uex}, and consider $\gamma = \tau - y_{12}$.
\begin{figure}[h]
\begin{minipage}{.3\linewidth}
\centering
\setlength{\unitlength}{0.254mm}
\begin{picture}(161,93)(20,-100)

        \allinethickness{0.254mm}\special{sh 0.3}\put(55,-30){\ellipse{4}{4}} 
        \allinethickness{0.254mm}\special{sh 0.3}\put(85,-30){\ellipse{4}{4}} 
        \allinethickness{0.254mm}\special{sh 0.3}\put(85,-45){\ellipse{4}{4}} 
        \allinethickness{0.254mm}\special{sh 0.3}\put(85,-60){\ellipse{4}{4}} 
        \put(20,-21){\textbf{\shortstack{$V'$}}} 
        \put(50,-21){\textbf{\shortstack{$P_1$}}} 
        \put(80,-21){\textbf{\shortstack{$P_2$}}} 
        \allinethickness{0.254mm}\special{sh 0.3}\put(85,-75){\ellipse{4}{4}} 
        \put(110,-21){\textbf{\shortstack{$P_3$}}} 
        \allinethickness{0.254mm}\special{sh 0.3}\put(25,-30){\ellipse{4}{4}} 
        \allinethickness{0.254mm}\special{sh 0.3}\put(55,-45){\ellipse{4}{4}} 
        \allinethickness{0.254mm}\special{sh 0.3}\put(55,-60){\ellipse{4}{4}} 
        \allinethickness{0.254mm}\special{sh 0.3}\put(55,-75){\ellipse{4}{4}} 
        \allinethickness{0.254mm}\special{sh 0.3}\put(115,-30){\ellipse{4}{4}} 
        \allinethickness{0.254mm}\special{sh 0.3}\put(115,-45){\ellipse{4}{4}} 
        \allinethickness{0.254mm}\special{sh 0.3}\put(115,-60){\ellipse{4}{4}} 
        \put(120,-66){\textbf{\shortstack{$y_{13}$}}} 
        \put(120,-51){\textbf{\shortstack{$y_7$}}} 
        \put(120,-36){\textbf{\shortstack{$y_4$}}} 
        \put(90,-66){\textbf{\shortstack{$y_9$}}} 
        \put(90,-51){\textbf{\shortstack{$y_6$}}} 
        \put(90,-36){\textbf{\shortstack{$y_3$}}} 
        \put(90,-81){\textbf{\shortstack{$y_{12}$}}} 
        \put(60,-96){\textbf{\shortstack{$y_{11}$}}} 
        \allinethickness{0.254mm}\special{sh 0.3}\put(55,-90){\ellipse{4}{4}} 
        \put(60,-81){\textbf{\shortstack{$y_{10}$}}} 
        \put(60,-66){\textbf{\shortstack{$y_8$}}} 
        \put(60,-51){\textbf{\shortstack{$y_5$}}} 
        \put(60,-36){\textbf{\shortstack{$y_2$}}} 
        \put(30,-36){\textbf{\shortstack{$y_1$}}} 
        \allinethickness{0.254mm}\path(110,-40)(110,-70)(80,-70)(80,-55)(50,-55)(50,-40)(20,-40)(20,-25)(80,-25)(80,-40)(110,-40) 
        \allinethickness{0.254mm}\path(50,-85)(80,-85)(80,-100)(50,-100)(50,-85) 
        \allinethickness{0.254mm}\path(110,-25)(110,-35)(110,-40)(140,-40)(140,-25)(110,-25) 
\end{picture}
\caption{$G$}
\end{minipage}%
\begin{minipage}{.3\linewidth}
\setlength{\unitlength}{0.254mm}
\begin{picture}(161,93)(20,-100)

        \allinethickness{0.254mm}\special{sh 0.3}\put(55,-30){\ellipse{4}{4}} 
        \allinethickness{0.254mm}\special{sh 0.3}\put(85,-30){\ellipse{4}{4}} 
        \allinethickness{0.254mm}\special{sh 0.3}\put(85,-45){\ellipse{4}{4}} 
        \allinethickness{0.254mm}\special{sh 0.3}\put(85,-60){\ellipse{4}{4}} 
        \put(20,-21){\textbf{\shortstack{$V'$}}} 
        \put(50,-21){\textbf{\shortstack{$P_1$}}} 
        \put(80,-21){\textbf{\shortstack{$P_2$}}} 
        \allinethickness{0.254mm}\special{sh 0.3}\put(85,-75){\ellipse{4}{4}} 
        \put(110,-21){\textbf{\shortstack{$P_3$}}} 
        \allinethickness{0.254mm}\special{sh 0.3}\put(25,-30){\ellipse{4}{4}} 
        \allinethickness{0.254mm}\special{sh 0.3}\put(55,-45){\ellipse{4}{4}} 
        \allinethickness{0.254mm}\special{sh 0.3}\put(55,-60){\ellipse{4}{4}} 
        \allinethickness{0.254mm}\special{sh 0.3}\put(55,-75){\ellipse{4}{4}} 
        \allinethickness{0.254mm}\special{sh 0.3}\put(115,-30){\ellipse{4}{4}} 
        \allinethickness{0.254mm}\special{sh 0.3}\put(115,-45){\ellipse{4}{4}} 
        \allinethickness{0.254mm}\special{sh 0.3}\put(115,-60){\ellipse{4}{4}} 
        \put(120,-66){\textbf{\shortstack{$y_{13}$}}} 
        \put(120,-51){\textbf{\shortstack{$y_7$}}} 
        \put(120,-36){\textbf{\shortstack{$y_4$}}} 
        \put(90,-66){\textbf{\shortstack{$y_9$}}} 
        \put(90,-51){\textbf{\shortstack{$y_6$}}} 
        \put(90,-36){\textbf{\shortstack{$y_3$}}} 
        \put(90,-81){\textbf{\shortstack{$y_{12}$}}} 
        \put(60,-96){\textbf{\shortstack{$y_{11}$}}} 
        \allinethickness{0.254mm}\special{sh 0.3}\put(55,-90){\ellipse{4}{4}} 
        \put(60,-81){\textbf{\shortstack{$y_{10}$}}} 
        \put(60,-66){\textbf{\shortstack{$y_8$}}} 
        \put(60,-51){\textbf{\shortstack{$y_5$}}} 
        \put(60,-36){\textbf{\shortstack{$y_2$}}} 
        \put(30,-36){\textbf{\shortstack{$y_1$}}} 
        \allinethickness{0.254mm}\path(50,-85)(80,-85)(80,-100)(50,-100)(50,-85) 
        \allinethickness{0.254mm}\path(20,-25)(135,-25)(135,-40)(105,-40)(105,-70)(80,-70)(80,-55)(50,-55)(50,-40)(20,-40)(20,-25) 
\end{picture}
\caption{$\tau'$}
\end{minipage}%
\begin{minipage}{.3\linewidth}
\setlength{\unitlength}{0.254mm}
\begin{picture}(161,93)(20,-100)

        \allinethickness{0.254mm}\special{sh 0.3}\put(55,-30){\ellipse{4}{4}} 
        \allinethickness{0.254mm}\special{sh 0.3}\put(85,-30){\ellipse{4}{4}} 
        \allinethickness{0.254mm}\special{sh 0.3}\put(85,-45){\ellipse{4}{4}} 
        \allinethickness{0.254mm}\special{sh 0.3}\put(85,-60){\ellipse{4}{4}} 
        \put(20,-21){\textbf{\shortstack{$V'$}}} 
        \put(50,-21){\textbf{\shortstack{$P_1$}}} 
        \put(80,-21){\textbf{\shortstack{$P_2$}}} 
        \allinethickness{0.254mm}\special{sh 0.3}\put(85,-75){\ellipse{4}{4}} 
        \put(110,-21){\textbf{\shortstack{$P_3$}}} 
        \allinethickness{0.254mm}\special{sh 0.3}\put(25,-30){\ellipse{4}{4}} 
        \allinethickness{0.254mm}\special{sh 0.3}\put(55,-45){\ellipse{4}{4}} 
        \allinethickness{0.254mm}\special{sh 0.3}\put(55,-60){\ellipse{4}{4}} 
        \allinethickness{0.254mm}\special{sh 0.3}\put(55,-75){\ellipse{4}{4}} 
        \allinethickness{0.254mm}\special{sh 0.3}\put(115,-30){\ellipse{4}{4}} 
        \allinethickness{0.254mm}\special{sh 0.3}\put(115,-45){\ellipse{4}{4}} 
        \allinethickness{0.254mm}\special{sh 0.3}\put(115,-60){\ellipse{4}{4}} 
        \put(120,-66){\textbf{\shortstack{$y_{13}$}}} 
        \put(120,-51){\textbf{\shortstack{$y_7$}}} 
        \put(120,-36){\textbf{\shortstack{$y_4$}}} 
        \put(90,-66){\textbf{\shortstack{$y_9$}}} 
        \put(90,-51){\textbf{\shortstack{$y_6$}}} 
        \put(90,-36){\textbf{\shortstack{$y_3$}}} 
        \put(90,-81){\textbf{\shortstack{$y_{12}$}}} 
        \put(60,-96){\textbf{\shortstack{$y_{11}$}}} 
        \allinethickness{0.254mm}\special{sh 0.3}\put(55,-90){\ellipse{4}{4}} 
        \put(60,-81){\textbf{\shortstack{$y_{10}$}}} 
        \put(60,-66){\textbf{\shortstack{$y_8$}}} 
        \put(60,-51){\textbf{\shortstack{$y_5$}}} 
        \put(60,-36){\textbf{\shortstack{$y_2$}}} 
        \put(30,-36){\textbf{\shortstack{$y_1$}}} 
        \allinethickness{0.254mm}\path(50,-85)(80,-85)(80,-100)(50,-100)(50,-85) 
        \allinethickness{0.254mm}\path(80,-55)(105,-55)(105,-70)(80,-70)(80,-55) 
\end{picture}
\caption{$R_O(\tau')$}
\end{minipage}
\end{figure}
\end{example}

If $y_{d+k} \notin U_O(\tau_i^k)$, let $\tau' = \gamma \cup s_O(\tau_i^k)$ (and recall that since $k>0$, we have seen that $s_O(\tau_i^k) < y_{d+k}$). Then $s_O(\tau') > s_O(\tau_i^k)$. Suppose $ y \in R_O(\tau')$ and $y\neq s_O(\tau_i^k)$. Again, $y \in \tau_i^k$. If $y > s_O(\tau')$, $y > s_O(\tau_i^k)$, so $y \in R_O(\tau_i^k)$. On the other hand, suppose $y \in U_O(\tau')$. If $y \in U_O(\tau_i^k)$, then $y \in R_O(\tau_i^k)$. Suppose $y \notin U_O(\tau_i^k)$. Then for some $q$, $\tau'$ contains all but one element, $b$, of $P_q$, $y \in P_q$, and $y >b$, but $\tau_i^k$ is missing at least 2 elements of $P_q$. As $s_O(\tau_i^k)$ is the only element of $\tau'$ not in $\tau_i^k$, $s_O(\tau_i^k)$ must be in $P_q$, and the only other element of $P_q$ not in $\tau_i^k$ must be $b$. In particular, $b \notin \tau_i^k$, so $b \geq s_O(\tau_i^k)$. But as $y \in U_O(\tau')$, $y > b$, so $y > s_O(\tau_i^k)$, and thus $y \in R_O(\tau_i^k)$. In particular, note that every element of $R_O(\tau') - s_O(\tau_i^k)$ is in $R_O(\tau_i^k)$, and as before we see that  $R_O(\tau') - s_O(\tau_i^k) \subseteq R_O(\tau_i^k) - y_{d+k}$. Hence $|R_O(\tau')| \leq |R_O(\tau_i^k)|$.

In either case, $\tau'$ is a facet of $skel_d(\Lambda)$ containing $\gamma$, and by construction must be equal to $\tau^r_t$ for some $r<k$. Since $|R_O(\tau')| \leq |R_O(\tau_i^k)|$, $\text{deg}(\sigma(\tau')) \leq \text{deg}(\sigma(\tau_i^k))$. If $r = 0$, $\sigma(\tau') = 1$, a divisor of $\sigma(\tau_i^k)$. Otherwise, $\sigma(\tau')$ is some monomial in $x_1, \ldots x_r$. Let $\mu$ be the reverse lexicographically largest divisor of $\sigma(\tau_i^k)$ whose degree is the same as that of $\sigma(\tau')$. Then $y_k$ divides $\mu$, and as the support of $\sigma(\tau')$ is in variables less than $x_k$, $\sigma(\tau') < \mu$. Thus (3) is proved.

\end{proof}


\section{Theorem \ref{novik}}
\label{generalize}
The proof of Theorem \ref{novik} is essentially that given by Novik in \cite{MR2122283} for the $p_i = p_j$ case, so we here give an abbreviated account with the necessary modifications, referring the reader to \cite{MR2122283} for full details.

Let $\Lambda = \Lambda(l; p_1, p_2, \ldots, p_m)$, and let $\Gamma$ be a $(d-1)$-dimensional Cohen-Macaulay subcomplex of $\Lambda$. Let $P_i$ and $V'$ be as defined in the previous section and label the vertices of $\Lambda$ with variables $x_1, x_2, \ldots, x_n$, ordered so that $x_i \in P_i$ for $1 \leq i \leq m$, and $x_i \in V'$ for $n-l+1 \leq i \leq
n$. Let $\mathbf{k}$ be a field and $\mathbf{k}[\mathbf{x}] = \mathbf{k}[x_1, x_2, \ldots , x_n]$. Recall that the Stanley-Reisner ideal of $\Gamma$, $I_{\Gamma}$, is the ideal generated by monomials $x_{i_1}x_{i_2}\ldots x_{i_s}$ such that $\{x_{i_1}, x_{i_2}, \ldots , x_{i_s} \}$ is \emph{not} a face of $\Gamma$.

For $g \in GL_n(\mathbf{k})$, $g$ defines an automorphism of $\mathbf{k}[\mathbf{x}]$ by $g(x_j) = \sum_{i=1}^n g_{ij}x_i$. We say $g$ possesses the Kind-Kleinschmidt condition if for every facet $\{x_{i_1}, x_{i_2}, \ldots , x_{i_r} \}$ of $\Gamma$, the submatrix of $g^{-1}$ obtained by taking the intersection of the rows numbered $i_1, i_2, \ldots , i_k$ with the last $d$ columns has rank $r$. For such a $g$, let $J(g) = gI_{\Gamma} + \gen{x_{n-d+1}, \ldots , x_n}$.

Finally, for $I$ an ideal in $\mathbf{k}[\mathbf{x}]$, let $Bs(g,I) = \{ \mu \in S(\infty^n) : \mu \notin \text{span}_{\mathbf{k}}(\{ \mu' : \mu \prec \mu' \} \cup I) \}$, where $\prec$ is the order given by $\mu \prec \mu'$ if either deg($\mu) < $ deg($\mu'$) or deg($\mu) = $ deg($\mu'$) and $\mu' < \mu$ in our original order on monomials (notice the reversal). The crux of the proof lies in the fact that $Bs(g,J(g))$ is a multicomplex, and that $F(Bs(g,J(g))) = h(\Gamma)$.  We additionally make use of the fact $Bs(g, J(g)) = Bs(g, gI_{\Gamma}) \cap S(\infty^{n-d})$. It thus suffices to construct a matrix $g$ satisfying the Kind-Kleinschmidt condition such that $Bs(g, gI_{\Gamma})$ does not contain $x_i^{p_i}$ for $1 \leq i \leq m$.

To do this we first pass to a larger field. Let $\mathbf{K} = \mathbf{k}(y_{ij},w_{ij},z_{ij})$ be the field of rational functions in $\sum_i (p_i -1) + l^2 + l(\sum_i p_i)$ variables, where $Y = (y_{ij})$, $W = (w_{ij})$ and $Z = (z_{ij})$ are $(\sum_i (p_i-1)) \times (\sum_i (p_i-1))$, $l \times l$, and $(\sum_i p_i) \times l$ matrices, respectively. Let $E = (E_{ij})$ be the $m \times (\sum_i (p_i-1))$ matrix where
\begin{equation} 
E_{ij} = 
   \begin{cases}     1 & \text{if } x_{j-m} \in P_i\\
      0 & \text{otherwise}.
   \end{cases}
\end{equation}
Define 
\begin{figure}[h!]
$g^{-1} = \left[ \begin{array}{cc}
               \left[ \begin{array}{cc}
                I_m & -EY\\
                0 & Y \end{array} \right] & Z\\
                0 & W \end{array} \right]$
so that $g = \left[ \begin{array}{cc}
               \left[ \begin{array}{cc}
                I_m & E\\
                0 & Y^{-1} \end{array} \right] & *\\
                0 & W^{-1} \end{array} \right]$.
\end{figure}

Now, for each $i$, $P_i \notin \Gamma$, so $I_{\Gamma}$ contains $\prod_{x_j \in P_i}x_j$. Then $gI_{\Gamma}$ contains 
\begin{equation*}
\prod_{x_j \in P_i}g(x_j) = \prod_{x_j \in P_i}\left(x_i + \sum_{k>m}g_{kj}x_k \right) = x_i^{p_i} + \sum \{ \alpha_{\mu}\mu : \mu \prec x_i^{p_i} \}.
\end{equation*}
Thus $x_i^{p_i} \notin Bs(g, gI_{\gamma})$, so $Bs(g, J) \subseteq S((p_1-1), (p_2-1), \ldots , (p_m-1), \infty^{n-d-m})$.

It remains only to show that $g$ satisfies the Kind-Kleinschmidt condition. Note that the $i^{th}$ row of $EY$ is equal to the sum of the rows of $Y$ indexed (in the larger matrix $g$) by $j > m$ such that $x_j \in P_i$. Since no facet of $\Lambda$ contains $P_i$, and the entries of $Y$, $W$ and $Z$ are algebraically independent, it then follows that for $\{x_{i_1}, x_{i_2}, \ldots , x_{i_d} \}$ a facet of $\Gamma$, the determinant of the submatrix of $g^{-1}$ defined by the intersection of the last $d$ columns and the rows numbered $i_1, \ldots i_d$ is non-zero, so the Kind-Kleinschmidt condition holds.

\section{Remarks}
\label{remarks}
 Note that the class of subcomplexes of $\Lambda(l;p_1, p_2, \ldots, p_m)$ is larger than that of complexes having proper $G$ action with corresponding orbit structure. Thus one does not expect our conditions to be sufficient for face numbers of complexes with proper group action. Indeed, in \cite{MR932113} Stanley showed necessary conditions on the $h$-vectors of centrally symmetric Cohen-Macaulay complexes not implied by our conditions, which were later generalized by Adin in \cite{MR1354672} to the case of Cohen-Macaulay complexes with proper $\Z_p$-action. It would be of interest to determine sufficient conditions in this more restricted case.

 Also, as mentioned in the introduction, there seems to be a close relationship between Corollary \ref{gen} and Theorem \ref{bfsgen}. In particular it may be possible to achieve a further generalization capturing both results as part of some larger phenomenon.

\bibliography{jb}
\end{document}